% 			Covering Letter
%
% AUTHOR: 	Robert H. Gilman
%		Department of Mathematics
%		Stevens Institute of Technology
%		Hoboken, NJ 07030
%		Office Phone: 201-216-5425
%		Department Secretary: 201-216-5448
%		Home Phone: 201-763-7058
%		E-mail: rgilman@vaxc.stevens-tech.edu
%
% TITLE: The Geometry of Cycles in the Cayley Diagram of a Group
%
% INTENDED PUBLICATION: Contemporary Mathematics
%			The Mathematical Legacy of Wilhelm Magnus
% 
% EDITOR: 	Kathryn Kuiken
%		Department of Mathematics
%		Polytechnic University
%		333 Jay Street
%		Brooklyn, NY 11201
%
% Note: All paperwork has been mailed to the editor.	
%
%
%---------------------------------------------------------------------
%	End of cover letter
%---------------------------------------------------------------------
%---------------------------------------------------------------------
% Beginning of AMS-TeX File
%---------------------------------------------------------------------
\input amstex
\documentstyle{amsppt}
%%%%%%%%%%%%%%%%%%%%%%%%%%%%%%%%% Start of private macros
\def\abs#1{\vert #1 \vert}
\def\i{^{-1}}
\def\c#1{\lceil #1 \rceil}
\def\G{\Gamma}
\def\part#1{\therosteritem#1}
\def\QED{$\square$}
%%%%%%%%%%%%%%%%%%%%%%%%%%%%%%%%% End of private macros

\NoBlackBoxes

\leftheadtext{ROBERT H. GILMAN}
\rightheadtext{THE GEOMETRY OF CYCLES IN THE CAYLEY DIAGRAM OF A GROUP}

\topmatter

\title The Geometry of Cycles in the Cayley Diagram of a Group\endtitle

\author Robert H. Gilman \endauthor

\address Department of Mathematics, Stevens Institute of Technology, Hoboken, New Jersey 07030 \endaddress

\email rgilman\@vaxc.stevens-tech.edu\endemail

\dedicatory Dedicated to the memory of Wilhelm Magnus \enddedicatory



\abstract A study of triangulations of cycles in the Cayley diagrams of finitely generated groups leads to a new geometric characterization of hyperbolic groups.
\endabstract

\thanks I thank the Institute for Advanced Study for its hospitality while part of this work was being done.\endgraf
This paper is in final form and no version of it will be submitted for publication elsewhere\endthanks

\endtopmatter

\document
\head 1. Introduction\endhead

During the past several years combinatorial group theory has received an infusion of ideas both from topology and from the theory of formal languages. The resulting interplay between groups, the geometry of their Cayley diagrams, and associated formal languages such as the language of all words defining the identity has led to several developments including the introduction of automatic groups \cite{Eps}, hyperbolic groups \cite{Gro}, and geometric and language-theoretic characterizations of finitely generated virtually free groups \cite{MS1} \cite{MS2}. A group is virtually free if it has a free subgroup of finite index and in particular if it is finite. In \cite{MS1} (together with \cite{Dun}) virtually free groups are characterized as those groups for which a finite set of diagonals suffices to triangulate all cycles in the Cayley diagram. Our goal here is to investigate triangulations of cycles for arbitrary finitely generated groups.

>From now on all groups under discussion are understood to be finitely generated and all sets of generators finite. $G$ will be a group, $S$ a set of generators closed under inverse, and $\Gamma$ the corresponding Cayley diagram. A path $\gamma$ of length $n \geq 1$ in $\Gamma$ is a sequence of group elements $g_0, \ldots ,g_{n}$ with an edge of $\Gamma$ from each $g_{i-1}$ to $g_i$. The label of $\gamma$ is the product in order of the labels of its edges. If $g_{n} = g_0$, then $\gamma$ is a cycle. A word in $S$ represents the identity in $G$ if and only if it is the label of a cycle. Finally $\c x$ stands for the least integer not less than $x$, and $\log x$ is to the base $2$.

\proclaim{Definition 1} A diagonal triangulation of a circle in the plane is obtained by distinguishing one or more points on the circle and joining them by chords in such a way that
\roster
\item No two chords meet in the interior of the circle;

\item The interior of the circle is divided into triangles;

\item Each arc of the circle between two neighboring distinguished points is one side of a triangle.
\endroster
Circles with one, two, or three distinguished points are considered to be triangulated without adding any chords. 
\endproclaim

\proclaim{Definition 2} A diagonal triangulation of a cycle $\gamma = g_0 \ldots g_n$ in $\Gamma$ is a diagonal triangulation of a circle with points $p_1, \ldots ,p_n$ distinguished and a corresponding labeling of these points, $g_1, \ldots ,g_n$, counterclockwise around the circle. For any chord $C$ with endpoints $p_i,p_j$, a word of minimum length representing $g_i\i g_j$ is called a label of $C$ in the direction from $p_i$ to $p_j$. The label of the arc of the circle from $p_{i}$ to $p_{j}$ is defined to be the label of the corresponding subpath  of $\gamma$. A triangulation of word $w$ in $S$ representing the identity in $G$ is a triangulation of any cycle with label $w$.
\endproclaim

>From now on we will say simply triangulation instead of diagonal triangulation. A triangulation of $\gamma$ makes the triangulated circle into a directed labeled graph in which arcs of the circle between adjacent distinguished points are directed counterclockwise, and each chord is construed as two associated edges, one in each direction. The label of any path in this directed graph from a point with label $g$ to one with label $g'$ is a word representing $g\i g'$.  

Recall that every choice of generators $S$ determines a metric on $G$ with distance $d(g,g')$ equal to the length of the shortest word in $S$ which represents $g\i g'$. This metric is extended to $\Gamma$ by making each edge isometric to the unit interval.

\proclaim{Definition 3} The length of a chord in a triangulation of a cycle is the distance in $G$ between the labels of its endpoints.  A $k$-triangulation is one in which all chords have length at most $k$, and $\Delta(n)$ is the minimum value of $k$ such that all cycles of length at most $n$ can be $k$-triangulated. In particular $\Delta(n) = 0$ for $1 \leq n \leq 3$. To display the dependence of $\Delta(n)$ on the generators $S$ write $\Delta_S(n)$.
\endproclaim

\medskip
Now let us consider triangulations of cycles in Cayley diagrams of arbitrary finitely generated groups. 
\proclaim{Theorem~A} 
For any group $G$ and set of generators of $G$, $\Delta(n) \leq \c{n/3}$. If for all sufficiently large $n$, $\Delta(n) < \c{n/3}$, then $G$ is finitely presented and satisfies an exponential isoperimetric inequality.
\endproclaim
\flushpar
The meaning of the last assertion of Theorem A is that there is a constant $c$ such that every word $w$ (in the generators of $G$) of length $n$ which defines the identity in $G$ is freely equivalent to the product of at most $c^n$ conjugates of the defining relators and their inverses. If this condition holds for one presentation of $G$, then it holds for all although the value of $c$ depends on the presentation. See \cite{Ger} for details.
\proclaim{Theorem~B} 
For any group $G$ the following conditions are equivalent.
\roster
\item For some set of generators and constant $K$, $\Delta (n) \le n/6 + K$;
\item G is a hyperbolic group;
\item For any set of generators there are constants $Q$ and $R$ such that $\Delta (n) \le Q\log (n) + R$.
\endroster
\endproclaim

\flushpar
The proof of this theorem relies on the characterization of hyperbolic groups by subquadratic isoperimetric inequalities \cite{Ol}, \cite{Pa}. By hyperbolic groups we mean the word hyperbolic groups of Gromov \cite{Gro}. $G$ is hyperbolic if there is a constant $\delta$ such that every geodesic triangle in $\Gamma$ (i.e., triangle whose sides are geodesic segments) has the property that each point on any one side is a distance at most $\delta$ from some point on one of the other two sides \cite{GH\rm, Proposition 21 of chapter 2}. 
The validity of this condition is independent of the choice of generators although the value of $\delta$ is not. Hyperbolic groups are finitely presented \cite{Sa\rm, Proposition 17}. Small cancellation groups satisfying the hypothesis $C'(1/6)$ or the hypotheses $C'(1/4)$ and $T(4)$ are hyperbolic \cite{Str}, but it is easy to see that $Z\times Z$ is not.

\proclaim{Theorem~C} If $G$ is hyperbolic, then either
\roster
\item $G$ is virtually free, and for any set of generators $\Delta(n)$ is bounded; or
\item $G$ is not virtually free, and for any set of generators there are constants $M, P, Q, R$ such that $M\log n +P \leq \Delta(n) \leq Q\log n + R$.
\endroster
\endproclaim

Our results show that roughly speaking $\Delta(n)$ is either linear, logarithmic, or bounded, and that the logarithmic case characterizes hyperbolic groups which are not virtually free. 

\head 2. Proof of Theorem~A and Theorem~B\endhead

Throughout this section $w=a_1\ldots a_n$, $n \ge 1$,  will stand for a word in the generators of $G$ representing the identity. For a fixed set of relators $\Cal R$, define $\alpha(w)$ to be the least integer such that $w$ is freely equivalent to a product of $\alpha(w)$ conjugates of relators in $\Cal R$. If there is no such product, $\alpha(w) = \infty$. For $n \ge 1$ define $\beta(n)=\max_{1\le \abs w \le n}\alpha(w)$.

\demo{Proof of Theorem A} To prove the first assertion of Theorem~A construct a triangulation of $w$ by picking distinguished points $p_1, \ldots ,p_n$ on a circle in the plane, labeling them with the group elements represented by the successive prefixes of $w$, and drawing chords 
\roster
\item From $p_n$ to $p_i$ for all $i$ with $2 \leq i \leq \c{n/3}$; and
\item From $p_{\c{n/3}}$ to $p_i$ for all $i$ with $\c{n/3}+2 \leq i \leq 2\c{n/3}$; and
\item From $p_{2\c{n/3}}$ to $p_i$ for all $i$ with $2\c{n/3}+2 \leq i \leq n$.
\endroster

For the second part of Theorem A assume $\Delta(n) < \c {n/3}$ for all $n \ge N \ge 4$, and let $\Cal R$ be the set of all relators of length at most $N$.  We will show $\beta(n) \leq 2^n$. Clearly $\beta(n) \leq 1$ if $n \leq N$; we may assume $n > N$. By induction on $n$, $\beta(n-1) \leq 2^{n-1}$, so we need only show $\alpha(w)\le 2^n$. By hypothesis $w$ has a triangulation with all chords of length less than $\c {n /3}$. Because chord length is an integer, all chords have length less than $n/3$. Any chord $C$ divides the circle into two arcs with labels $w_3w_1$ and $w_2$ where $w=w_1w_2w_3$. 
Choosing a label $v$ for $C$ in the appropriate direction, we see that $w_1vw_3$ and $w_2v\i$ are both words representing the identity in $G$. It follows that $\alpha(w) \le \alpha(w_1vw_3)+\alpha(w_2v\i )$. Since labels of chords have minimum length, $\abs v \leq \abs {w_1w_3}$ and $\abs v \leq \abs {w_2}$. If both these inequalities are strict, then by induction $\alpha(w) \leq \alpha(w_1vw_3) + \alpha(w_2v\i) \leq 2^{n-1} + 2^{n-1} \leq 2^n$.

It remains to find a chord $C$ which divides the circle into two arcs both longer than $C$. Since $\abs w > N \ge 4$, the triangulation does have chords. Any chord divides the circle into two arcs of length, say, $d$ and $e$. Pick $C$ with $d$ as large as possible subject to $d \le e$. We claim $d \ge n/3$. To see this observe that $C$ is one side of a triangle $T$ which has its third vertex on the arc of length $e$. This third vertex divides that arc into two shorter arcs of length $d'$ and $d''$ with $e=d'+d''$, and our claim will follow from $d'\le d$ and $d'' \le d$. 
If $d' =1$, then clearly $d'\le d$. Otherwise the corresponding side of $T$ is not a subarc of the circle but a chord $C'$ which divides the circle into arcs of length $d'$ and $d+d''$. Since $d+d'' > d$, our choice of $C$ implies first that $d+d'' > d'$ and consequently that $d'\le d$. Thus $d' \le d$ in all cases; and as $d'' \le d$ by symmetry, our claim is valid. Since $C$ has length less than $n/3$, it is the desired chord.
\QED \enddemo

\demo{Proof of Theorem B} The proof that \part3 implies \part1 is straightforward and is omitted. Assume \part1 holds; that is, $\Delta (n) \le n/6 + K$. By \cite{Ol} or \cite{Pa} conclusion \part2 holds once we know that $G$ has a subquadratic isoperimetric inequality, i.e., $\lim_{n \to \infty}\beta(n)/n^2 = 0$. Take $N$ to be an integer larger than $\max\{3, 1000K\}$. We will show that $\beta(n) \le \beta(N)n^{1.9}$.

If $n \le N$, there is nothing to prove, so assume $n >N$; by induction on $n$ it suffices to show $\alpha(w)\leq \beta(N)n^{1.9}$. As in the proof of Theorem A find a chord $C$ which divides the circle into two arcs of lengths $e \ge d \ge n/3$. Note that $e=n-d\le 2n/3$ and $d \le n/2$. Since $C$ has length at most $\Delta(n)$, $\alpha(w) \le \beta(d+\Delta(n))+\beta(e+\Delta(n)) \le \beta(d + n/6 + K)+\beta(e+ n/6 + K)$. Because $d+n/6+K \le e + n/6 +K \le 2n/3 + n/6 + K < n$, induction on $n$ yields
$$\align 
\alpha(w) &\le \beta(d+n/6+K) + \beta(n - d+n/6+K) \\
	&\le \beta(N)n^{1.9}\bigl ( (d/n + 1/6 + .001)^{1.9} + (1 - d/n + 1/6 + .001)^{1.9}\bigr ) \\
	&\le \beta(N)n^{1.9}\cr
\endalign$$
where $ (d/n + 1/6 + .001)^{1.9} + (1 - d/n + 1/6 + .001)^{1.9} \le 1$ follows from $1/3 \le d/n \le 1/2$. 

\medskip
To complete the proof of Theorem B we will show that \part2 implies \part3 by proving that if $G$ is hyperbolic, then $\Delta(n) \le C\log(n)$ for some constant $C$. By assumption $G$ is $\delta/4$ hyperbolic for some $\delta$ (this odd choice of $\delta$ is made to correspond to the hypothesis of \cite{CDP\rm, Lemma~1.6 of Chapter 3}, which will be employed later). Since it does no harm to increase $\delta$, assume $\delta \ge 1$, and take $C > 10\delta$. 

Let $\gamma= g_0\ldots g_n$ be a cycle in $\Gamma$. As $\Delta(n) = 0$ for $n=1,2,3$, we may assume $n\ge 4$. Choose distinguished points $p_1 \ldots , p_n$ on a circle as in Definition 1, and give each $p_i$ the label $g_i$ as in Definition 2. Start constructing a triangulation by adding a chord from $p_n$ to $p_2$. As this chord has length at most 2, we are done if $n=4$. Otherwise it suffices to show that whenever a chord of length at most $C\log(n)$ has endpoints $p_i, p_j$ with $3 \le j-i$, then we can add a chord from $p_i$ to $p_{j-1}$ or from $p_{i+1}$ to $p_{j}$ or we can add chords from $p_i$ and $p_j$ to $p_k$ for some $k$ with $i+2 \le k \le j-2$. In other words 
\roster
\item $d(g_i,g_{j-1})\le C\log(n)$; or
\item $d(g_{i+1},g_j)\le C\log(n)$; or
\item $d(g_i,g_k), d(g_j,g_k) \le C\log(n)$ for some $k$ with $i+2 \le k \le j-2$.
\endroster 

Suppose \part1 and \part2 do not hold. Thus $d(g_i, g_j) > C\log(n)-1$. Let $\gamma'$ be a geodesic in $\Gamma$ from $g_i$ to $g_j$ and consider a ball of radius $r=(C/2)\log(n) -2$ around the midpoint $x$ of $\gamma'$. If the part of $\gamma$ from $g_i$ to $g_j$ intersects the ball, then there is a vertex $g_k$ on $\gamma$ with $i \le k \le j$ and with $g_k$ a distance at most $(C/2)\log(n) + r +1 \le C\log(n) - 1$ from each endpoint. Since $d(g_i,g_j) > C\log (n) -1$, we have $i < k < j$; and one of the conditions above must hold.

If $\gamma$ does not intersect the ball, then by \cite{CDP\rm, Lemma~1.6 of Chapter 3}, $j-i \geq \delta(2^{(r/\delta) -1} - 2)$. As $\delta \ge 1$ and $i \ge 2$, we have $n \ge j \ge j-i+2 \geq 2^{(r/\delta) -1}= 2^{C\log(n)/2\delta - 2/\delta-1} \ge n^{C/2\delta}2^{-3} \ge n^{5}/8$. But $n > n^5/8$ is impossible as $n \ge 2$. \QED \enddemo

\head 3. Proof of Theorem~C
\endhead
\proclaim{Lemma~3.1} 
If $S$ and $S'$ are generating sets for $G$, there is a constant $K$ such that $\Delta_S(n) \le K\Delta_{S'}(Kn) + K$. 
\endproclaim
\demo{Proof} First note that given a triangulation of the circle in the sense of Definition 1,  with distinguished points $p_1\ldots p_n$, $n > 3$, we can generate a new triangulation with $n-1$ distinguished points by allowing a point $p_i$ to move counterclockwise on the circle until it becomes identified with $p_{i+1}$. Here $i+1,i-1$ etc.\ are understood modulo $n$. Let $p_j$ be the third vertex of the triangle whose other two vertices are $p_i,p_{i+1}$. 
\roster
\item If $j \notin \{i+2,i-1\}$, then the chord $\overline{p_ip_j}$ is identified with the chord $\overline{p_{i+1}p_j}$, and the triangle $\overline{p_ip_{i+1}p_j}$ disappears. All other triangles with $p_i$ as a vertex have that vertex replaced by $p_{i+1}$.
\item If $j=i+2$, then the chord $\overline{p_ip_{i+2}}$ is identified with the edge $\overline{p_{i+1}p_{i+2}}$, and the triangle $\overline{p_ip_{i+1}p_{i+2}}$ disappears. All other triangles with $p_i$ as a vertex have that vertex replaced by $p_{i+1}$.
\item If $j=i-1$, then the original triangulation has the triangle $\overline{p_{i-1}p_ip_{i+1}}$; and the new triangulation is obtained by removing the chord $\overline{p_{i-1}p_{i+1}}$. 
\endroster

Now choose $K_1$ so that each generator in $S$ can be expressed as a word of length at most $K_1$ in $S'$ and vice-versa. Let $\Gamma$ and $\Gamma'$ be the Cayley diagrams of $G$ with respect to $S$ and $S'$ respectively, and take $d$ and $d'$ be the corresponding metrics. We have $(1/K_1) d'(g,h) \le d(g,h) \le K_1 d'(g,h)$. Suppose $\gamma = g_0 \ldots g_n$ is a cycle of length $n$ in $\Gamma$. If $n \le 3$, then there is nothing to prove as $\Delta_S(n)=0$. Thus we may assume $n > 3$. Since $d(g_{i-1},g_{i}) \leq 1$ implies $d'(g_{i-1},g_{i}) \leq K_1$, $\gamma$ can be expanded to a cycle $\gamma'$ of length $K_1n$ or less in $\Gamma'$ by interpolating at most $K_1-1$ additional vertices between each $g_{i-1}$ and $g_{i}$. 

Consider a $\Delta_{S'}(K_1n)$-triangulation of $\gamma'$ with distinguished points $p_1\ldots p_n$ corresponding to the original cycle $\gamma$ and additional points $q_{i,1}\ldots q_{i,j(i)}$, $j(i)<K_1$, corresponding to the additional vertices between $g_i$ and $g_{i+1}$. Modify this triangulation as above so that the points $q_{i,j}$ are all identified with $p_{i+1}$ to obtain a triangulation with distinguished points $p_1, \ldots , p_n$. As each point $q_{i,j}$ moves a distance at most $K_1$ on the circle, each chord lengthens by at most $2K_1$ in the metric $d'$. We obtain a $\Delta_{S'}(K_1n)+2K_1$-triangulation in terms of $d'$ and consequently a $K_1\Delta_{S}(K_1n)+2K_1^2$-triangulation of $\gamma$ with respect to $d$.
\QED \enddemo
\proclaim{Lemma~3.2} Suppose $G$ is a free product of $H$ and $K$ with finite subgroups of $H$ and $K$ amalgamated or $G$ is an HNN extension with base $H$ and two finite subgroups of $H$ associated. There is a set of generators for $H$ which extends to a set of generators for $G$ in such a way that the subgroup $H$ is isometrically embedded in $G$ with respect to the corresponding metrics $d_H, d_G$.  
\endproclaim
\demo{Proof} We treat the free product case first. Choose a set of generators $S_H$ for $H$ which includes every element of its amalgamated subgroup (and is closed under inverse). Choose $S_K$ likewise for $K$. $S_G=S_H \cup S_K$ is a set of generators for $G$. Clearly $d_H(h_1,h_2) \ge d_G(h_1, h_2)$. To prove the reverse inequality suppose $u$ and $w$ are words of minimum length in $S_H$ and $S_G$ respectively both representing $h=h_1\i h_2$. The word $w$ factors uniquely as $w=w_1\ldots w_n$ in such a way that $w_i$ and $w_{i+1}$ are words in different alphabets from $\{S_H, S_K\}$. Among all words of minimum length representing $h$ we choose $w$ with $n$ minimum. It suffices to show that 
$$
\abs u \le \abs w.\tag 3.1
$$
If $n>1$, then no subword $w_i$ represents an element of an amalgamated subgroup. Otherwise $w_i$ could be replaced by a single generator from the other alphabet to obtain a new word which represented the same element of $h$, was no longer than $w$, and whose factorization as a product of words in different alphabets had fewer than $n$ terms. The word $u\i w_1\ldots w_n$ represents the identity, and the first term of its factorization into a product of words from different alphabets is either $u\i$ or $u\i w_1$ depending on whether $w_1$ is a word in $S_K$ or $S_H$. 
By the normal form theorem for free products with amalgamation \cite{LS\rm, Chapter 4} one of the terms in the factorization must represent an element of an amalgamated subgroup, and the only possibility is $u\i$ or $u\i w_1$ respectively. In the first case $\abs u \le 1$ because it represents a word in the amalgamated subgroup of $H$, and (3.1) follows directly. Likewise in the second case the group element represented by $u\i w_1$ is also represented by a word of length at most $1$ whence $\abs u \le \abs {w_1} + 1$. Consequently (3.1) holds unless $w=w_1$; but then as $w_1$ is a word in $S_H$, (3.1) holds by choice of $u$.    

The HNN case is similar to the free product case. Suppose that is $G$ an HNN extension with base $H$ and stable letter $t$ and that $A$ and $B$ are the associated subgroups of $H$ with $t\i A t = B$. Choose generators $S_H$ for $H$ with $A \cup B \subset S_H$, and let $S_G=S_H\cup \{t,t\i\}$. It suffices to prove (3.1) when $u$ and $w$ are chosen as before. 
More precisely any word $w$ in $S_G$ factors uniquely as $w=w_0t^{\epsilon_1} \ldots t^{\epsilon_n}w_n$ where $\epsilon_i=\pm1$, the $w_i$'s are words (possibly empty words) in $S_H$; and $w$ is chosen with $n$ minimum among words of minimum length representing $h$. It is straightforward to check that this factorization does not include any subsequences of the form $t\i w_i t$ with $w_i$ representing an element of $A$ or $t w_i t\i$ with $w_i$ representing an element of $B$. Applying Britton's Lemma \cite{LS\rm, Chapter 4} to $(u\i w_0) t^{\epsilon_1} \ldots t^{\epsilon_n}w_n$, we conclude that $n=0$ whence $\abs u \le \abs{w_0}= \abs w$ by choice of $u$.
\QED \enddemo

\demo{Proof of Theorem C} If $G$ is virtually free, then conclusion \part1  holds by \cite{MS1}; the upper bound of \part2) comes from Theorem B.  Thus it suffices to show 
$$
M\log n +P \leq \Delta(n)\tag 3.2
$$
when $G$ is hyperbolic but not virtually free. First we show that although the constants $M$ and $P$ may change from one generating set to another, the validity of (3.2) is independent of the choice of generators for $G$. 

Suppose that (3.2) holds for a set of generators $S$. If $S'$ is another set, then  $M\log n +P \leq \Delta_{S}(n) \le K\Delta_{S'}(Kn) + K$ by Lemma 3.1. We claim $M'\log m +P' \leq \Delta_{S'}(m)$ for some constants $M',P'$. If this inequality holds for all but finitely many $m$, then with a change of $P'$ it holds for all $m$. Thus we may assume $m \ge 2K$ and choose $n\ge 2$ so that $nK \le m < (n+1)K$. We obtain $m \le Kn+K$ and $M\log n +P \leq  K\Delta_{S'}(Kn) + K\le K\Delta_{S'}(m) + K$. The desired inequality follows in a straightforward way. 

By the preceding argument it suffices to show that (3.2) holds for one set of generators. Since hyperbolic groups are finitely presented, we may use induction on the accessibility length of $G$ \cite{Dun}. The accessibility length of $G$ is the length of the longest series $G=G_0 \supset G_1 \supset \ldots \supset G_n$ such that each $G_i$ has a decomposition as a nontrivial free product with amalgamation or as an HNN extension where one of the factors or the base is $G_{i+1}$ and the amalgamated or associated subgroups are finite. 
A free product with amalgamation is nontrivial if the amalgamated subgroups are both proper. As $G$ is not virtually free, the results of \cite{Sta} imply that the number of ends of $G$ is either 1 or $\infty$, and in the latter case $G$ has a decomposition of the type mentioned. 

Suppose $G$ is a nontrivial free product $G=H*_FK$ with $F$ finite. Choose the generators of Lemma 3.2. It is an immediate consequence of that lemma that $H$ and $K$ are hyperbolic. If $H$ and $K$ are virtually free, then by \cite{Gre} or \cite{KPS} so is $G$. Thus we may assume $H$ is not virtually free. As $H$ has shorter accessibility length than $G$, the induction assumption yields $\Delta_H(n) \ge M\log n + P$, 
and $\Delta_G(n) \ge M\log n + P$ follows from Lemma 3.2. Likewise if $G$ is an HNN extension with base $H$ and finite associated subgroups, then as before $H$ is hyperbolic but not virtually free by \cite{Gre} or \cite{KPS}. The induction hypothesis yields $\Delta_H(n) \ge M\log n + P$, and the desired inequality follows.

It remains to deal with the case in which $G$ has one end. Let $\G$ be the Cayley diagram of $G$ with respect to some set of generators $S$, and let $d$ be the corresponding word metric. Since $G$ is hyperbolic, geodesic triangles in $\G$ are $\delta$-thin for some $\delta $ \cite{GH\rm, Definition 16}. The meaning of $\delta$-thin is that there is a map from the perimeter of the triangle to three lines in the Euclidean plane with a common endpoint such that 
\roster
\item The vertices of the triangle are mapped onto the other endpoints of the lines;
\item The restriction of $f$ to each side of the triangle is an isometry;
\item Points with the same image under $f$ are a distance at most $\delta$ apart.
\endroster
Since increasing $\delta$ does no harm, take $\delta$ to be a positive integer. Pick a geodesic segment $\gamma$ in $\G$ of length $2n$ with ends $g,g'$ and midpoint $1$. As $G$ has one end, there is a path $\gamma'=g_1, g_2, \ldots g_N$ in $\G$ such that $g_1=g$, $g_N=g'$, and  $d(1, g_i)\ge n$ for all $i$. 

We will find a new path $\gamma''$ from $g$ to $g'$. For each $i$, $1 < i < N$, pick a geodesic segment from $1$ to $g_i$. Let $h_i$ be the vertex on this segment with $d(1,h_i)=n$ and $\gamma_i$ the subsegment from $1$ to $h_i$. Define $h_1=g_1$ and $h_N=g_N$ and take $\gamma_1$ and $\gamma_N$ be the subsegments of $\gamma$ from $1$ to $g_1$ and $g_N$ respectively. For each $i$, $1 \leq i < N$, consider the geodesic triangle with vertices $1,g_i, g_{i+1}$, and whose edges are the geodesic segments from $1$ to $g_i$ and $g_{i+1}$ previously chosen together with the edge in $\gamma$ from $g_i$ to $g_{i+1}$.
As this triangle is $\delta$-thin and the side opposite vertex $1$ has length $1$, it follows in a straightforward way that $d(h_i,h_{i+1}) \le 2\delta + 1$. Construct a path from $h_1$ to $h_N$ by joining each $h_i$ to $h_{i+1}$ with a geodesic segment of length at most $2\delta + 1$. Clearly the distance from $1$ to any point on this path is at least $n-(2\delta+1)$. 
As the labels of the $\gamma_i$'s are words of length $n$ in $S$, there are at most $\abs{S}^n$ distinct $\gamma_i$'s and hence at most that many distinct $h_i$'s. Thus by deleting loops from the path just constructed, we obtain a path $\gamma''$ from $h_1$ to $h_N$ of length at most $(2\delta + 1)\abs{S}^n$. Further the distance from $1$ to any  point on $\gamma''$ is at least $n-(2\delta+1)$. 
 
Consider any triangulation of the cycle formed by $\gamma$ and $\gamma''$. This cycle has length at most 
$$
L(n)=2n + (2\delta + 1)\abs{S}^n\tag 3.3
$$ 
By \cite{MS\rm, Lemma 5} any diagonal triangulation has the property that if the circle is divided into three arcs each beginning and ending at distinguished points, then there is a triangle with vertices on each arc. Thus there is a triangle with vertex $h$ on $\gamma''$, and vertices on $\gamma_1$ and $\gamma_N$. It follows that $d(1,h)$ is at most equal to the sum of the lengths of two sides of this triangle whence $n-2\delta-1 \le d(h,1) \le 2\Delta (L(n))$.    

Pick any $m \ge 2+(2\delta + 1)\abs{S}$ and choose $n$ with $L(n)\le m \le L(n+1)$. From the preceding paragraph $n/2-\delta -1/2 \le \Delta(L(n)) \le \Delta(m)$. On the other hand as $m \le L(n+1)$, (3.3) yields $M\log m + P \le n/2-\delta -1/2 \le \Delta(m)$ for some constants $M,P$. 
\QED \enddemo

\head Acknowledgement \endhead
I thank Panagiotis Papasoglu for showing me how to improve an earlier version of Theorem C.

\enddocument